\documentclass[12pt]{article}
\usepackage{amsmath,amsthm,amssymb}
\begin{document}
\newtheorem{thm}{Theorem}
\newtheorem{pro}[thm]{Proposition}
\newtheorem{cor}[thm]{Corollary}
\newtheorem{lem}[thm]{Lemma}
\newtheorem{df}[thm]{Definition}
\newtheorem{rem}[thm]{Remark}
\newcommand{\End}{\mathop{\mathrm{End}}\nolimits}
\newcommand{\Id}{\mathop{\mathrm{Id}}\nolimits}
\newcommand{\Ric}{\mathop{\mathrm{Ric}}\nolimits}
\title{Integral Geometry and Hamiltonian volume minimizing property of
a totally geodesic Lagrangian torus in $S^{2} \times S^{2}$}
\author{Hiroshi Iriyeh\footnote{The first author is supported by Grant-in-Aid
for JSPS Fellows 08889.}, Hajime Ono\footnote{The second author is supported
by Grant-in-Aid for JSPS Fellows 08832.} and Takashi Sakai}
\date{}
\maketitle

\begin{abstract}
We prove that the product of equators $S^{1} \times S^{1}$ in
$S^{2} \times S^{2}$ is globally volume minimizing under Hamiltonian
deformations.

\end{abstract}

\section{Introduction and main results}

In 1990, Y.\ G.\ Oh \cite{Oh} introduced the notion of {\it global Hamiltonian
stability} of minimal Lagrangian submanifolds in a K\"{a}hler manifold and
posed the following conjecture:\newline
{\bf Conjecture (Oh).}\ Let $M$ be a K\"{a}hler-Einstein manifold with an
involutive anti-holomorphic isometry $\tau$.
Suppose that the fixed point set of $\tau$
\begin{eqnarray*}
 L:=\mbox{Fix}\ \tau
\end{eqnarray*}
is also a compact Einstein manifold with positive Ricci curvature.
Then for any Hamiltonian isotopy $\rho \in \mbox{Ham}(M)$ of $M$,
we have
\begin{eqnarray*}
 \mbox{vol}(\rho(L)) \geq \mbox{vol}(L).
\end{eqnarray*}

Kleiner and Oh \cite{Oh} proved that this conjecture is true for the case
${\mathbb R}P^{n} \subset {\mathbb C}P^{n}$ (see also \cite{Gi}).
\begin{thm}[Kleiner-Oh]
The standard ${\mathbb R}P^{n} \subset {\mathbb C}P^{n}$ has the least volume
among all its images under Hamiltonian isotopies.
\end{thm}
This is the only known example such that the conjecture has been proved
affirmatively.

Important examples of K\"{a}hler-Einstein manifolds admitting an involutive
anti-holomorphic isometry are Hermitian symmetric spaces.
Let $M$ be a Hermitian symmetric space of compact type and $\tau$ be a
canonical involution on $M$.
Then
\begin{eqnarray*}
 L:=\mbox{Fix}\ \tau
\end{eqnarray*}
is a totally geodesic Lagrangian submanifold in $M$ (which is called a
{\it real form} of $M$).
It is interesting to verify the conjecture for such a pair $(M,L)$.

In this paper, we shall prove that the same statement as the conjecture is
true in the case of $(S^2 \times S^2 \cong Q_{2}({\mathbb C}),S^1 \times S^1)$
although the Lagrangian surface $S^1 \times S^1$ is {\it flat}.
More precisely,
\begin{thm}
Let $L:=S^1 \times S^1$ be a totally geodesic Lagrangian torus in
$(S^2 \times S^2,\omega_{0} \oplus \omega_{0})$, where $\omega_{0}$ denotes
the standard K\"{a}hler form of $S^2(1) \cong {\mathbb C}P^1$.
Then for any Hamiltonian isotopy $\rho \in \mathrm{Ham}(S^2 \times S^2)$,
we have
 \begin{eqnarray*}
  \mathrm{vol}(\rho(L)) \geq \mathrm{vol}(L).
 \end{eqnarray*}
\end{thm}

Our proof is based on the following Lagrangian intersection theorem
(\cite{Oh3},\cite{Oh1} and \cite{Oh2}) and a new Poincar\'e formula for
Lagrangian surfaces in $S^2 \times S^2$.
\begin{thm}[Oh]
Let $(M,\omega)$ be a compact symplectic manifold such that there exists an
integrable almost complex structure $J$ for which the triple $(M,\omega,J)$
becomes a compact Hermitian symmetric space.
Let $L=\mathrm{Fix}\ \tau$ be the fixed point set of an anti-holomorphic
involution isometry $\tau$ on $M$.
Assume that the minimal Maslov number of $L$ is greater than or equal to $2$.
Then for any Hamiltonian isotopy $\rho \in \mathrm{Ham}(M)$ of $M$ such that
$L$ and $\rho(L)$ intersect transversally, the inequality
 \begin{eqnarray}
  \sharp (L \cap \rho (L)) \geq
   \sum_{i=0}^{\dim L} \mathrm{rank} H_i(L, {\mathbb Z}/2{\mathbb Z})
 \end{eqnarray}
holds.
\end{thm}
Since the minimal Maslov number of $S^1 \times S^1 \subset S^2 \times S^2$
is $2$, the assumption of the above theorem is satisfied in our case.
\begin{pro}
Let $N$ and $L$ be surfaces of $S^2 \times S^2$.
Suppose that $N$ is Lagrangian and $L$ is a product of curves in $S^2$.
Then the following inequality holds$:$
 \begin{eqnarray}
  4 \pi \mathrm{vol}(N) \mathrm{vol}(L)
  \leq \int_{SO(3) \times SO(3)} \sharp (N \cap gL)d\mu(g)
  \leq 16 \mathrm{vol}(N) \mathrm{vol}(L).
 \end{eqnarray}
\end{pro}

This formula is interesting in its own right.
We remark the equality condition of the inequality $(2)$.
The first equality of $(2)$ is fulfilled by, for example, a Lagrangian
embedding $S^2 \ni z \mapsto (z,-z) \in S^2 \times S^2$.
The second equality of $(2)$ holds if and only if the Lagrangian surface $N$
is also a product of closed curves in $S^2$.

\section{Poincar\'e formula in Riemannian homogeneous spaces}

Here we shall review the generalized Poincar\'e formula
in Riemannian homogeneous spaces obtained by Howard \cite{Howard}.

Let $U$ be a finite dimensional real vector space with an inner product,
and $V$ and $W$ vector subspaces of dimension $p$ and $q$ in $U$,
respectively.
Take orthonormal bases $v_1, \ldots, v_p$ and $w_1, \ldots, w_q$
of $V$ and $W$, and define
$$
\sigma(V,W)
= \| v_1 \wedge \cdots \wedge v_p \wedge w_1 \wedge \cdots \wedge w_q \|,
$$
which is the angle between $V$ and $W$.

Let $G$ be a Lie group and $K$ a closed subgroup of $G$.
We assume that $G$ has a left invariant Riemannian metric
which is also invariant under elements of $K$.
This metric induces a $G$-invariant Riemannian metric on $G/K$.
For $x$ and $y$ in $G/K$ and vector subspaces $V$ in $T_x(G/K)$ and
$W$ in $T_y(G/K)$,
we define $\sigma_K^{}(V,W)$, the angle between $V$ and $W$, by
$$
\sigma_K(V,W)
= \int_K \sigma( (dg_x)_o^{-1} V, dk_o^{-1} (dg_y)_o^{-1} W ) d\mu_K^{}(k)
$$
where $g_x$ and $g_y$ are elements of $G$ such that $g_xo=x$ and $g_yo=y$.
Here we denote by $o$ the origin of $G/K$.
\begin{thm}[Howard]
Let $G/K$ be a Riemannian homogeneous space and assume that $G$ is unimodular.
Let $N$ and $L$ be submanifolds of $G/K$
with $\dim N + \dim L \geq \dim (G/K)$.
Then
$$
\int_G {\rm vol}(N \cap gL) d\mu_G^{}(g)
= \int_{N \times L} \sigma_K(T_x^\perp N, T_y^\perp L) d\mu(x,y)
$$
holds.
\end{thm}

The linear isotropy representation induces an action of $K$
on the Grassmannian manifold $G_p(T_o(G/K))$
consisting of all $p$ dimensional subspaces
in the tangent space $T_o(G/K)$ at $o$ in a natural way.
Although $\sigma_K(T_x^\perp N, T_y^\perp L)$ is defined as an integral on $K$,
we can consider that it is defined as an integral on an orbit of $K$-action on
the Grassmannian manifold.
So $\sigma_K(\cdot, \cdot)$ can be regarded as a function
defined on the product of the orbit spaces of such $K$-actions.
In the case where $G/K$ is a real space form,
$\sigma_K(T_x^\perp N, T_y^\perp L)$ is constant
since $K$ acts transitively on the Grassmannian manifold.
This implies that the Poincar\'e formula is expressed
as a constant times of the product of the volumes of $N$ and $L$.
In general, such $K$-actions are not transitive.
However, if we can define an invariant for orbits of this action,
which is called an {\em isotropy invariant},
then using this we can express the Poincar\'e formula more explicitly.
>From this point of view,
Tasaki \cite{Tasaki} introduced the multiple K\"ahler angle,
which is the invariant for the actions of unitary groups.

\section{Poincar\'e formula for Lagrangian surfaces in $S^{2} \times S^{2}$}

In this section we define isotropy invariants
for surfaces in $S^2 \times S^2$,
and give a concrete expression of the Poincar\'e formula
for its Lagrangian surfaces.

Let $G$ be the identity component of the isometry group of $S^2 \times S^2$,
that is, $G=SO(3) \times SO(3)$.
Then the isotropy group $K$ at $o = (p_1, p_2)$ in $S^2 \times S^2$
is isomorphic to $SO(2) \times SO(2)$,
and $S^2 \times S^2$ is expressed as a coset space $G/K$.
Assume thet $G$ is equipped with an invariant metric normalized so that
$G/K$ becomes isometric to the product of unit spheres.
We decompose the tangent space $T_o(G/K)$ as
$$
T_o(G/K) = T_{p_1}(S^2) \oplus T_{p_2}(S^2).
$$
We take orthonormal bases $\{ e_1, e_2 \}$ and $\{ e_3, e_4 \}$
of $T_{p_1}(S^2)$ and $T_{p_2}(S^2)$, respectively,
then a complex structure on $T_o(G/K)$ is given by
$$
Je_1 = e_2,\ Je_2 = -e_1,\ Je_3 = e_4,\ Je_4 = -e_3.
$$

We consider the oriented $2$-plane Grassmannian manifold
$\tilde G_2(T_o(G/K))$.
Take an origin $V_o := \mbox{span} \{ e_1, e_2 \}$,
and express $\tilde G_2(T_o(G/K))$ as a coset space
$$
\tilde G_2(T_o(G/K)) = SO(4) / (SO(2) \times SO(2)) =: G'/K'.
$$
Now we study the $K$-action on $\tilde G_2(T_o(G/K))$,
and define isotropy invariants.
In this case the actions of $K$ and $K'$ on $\tilde G_2(T_o(G/K))$
are equivalent by $\mbox{Ad} : K \to K'$.
Therefore it suffices to consider the orbit space
of the isotropy action of $\tilde G_2(T_o(G/K))$.
It is well known that the orbit space of the isotropy action
of a symmetric space of compact type can be identified
with a fundamental cell of a maximal torus.
Hence we can define the isotropy invariant by a coordinate of a maximal torus.
We denote by $\frak g'$ and $\frak k'$ the Lie algebra of $G'$ and $K'$,
respectively.
Then we have a canonical orthogonal direct sum decomposition
$\frak g' = \frak k' \oplus \frak m'$ of $\frak g'$, where
$$
\frak m' = \left\{ \left( \begin{array}{cc} O & X \\ -^tX & O \end{array} 
\right) \ \bigg| \ X \in M_2(\mathbb R) \right \}.
$$
We take a maximal abelian subspace $\frak a'$ of $\frak m'$ as follows:
$$
\frak a' = \left\{
\left( \begin{array}{cc} O & X \\ -^tX & O \end{array} \right) \ \bigg| \
X = \left( \begin{array}{cc} \theta_1 & 0 \\ 0 & \theta_2 \end{array} \right),
\theta_1, \theta_2 \in \mathbb R \right \}.
$$
Then the set of positive restricted roots of $(\frak g', \frak k')$
with respect to $\frak a'$ is
$$
\Delta = \{ \theta_1 + \theta_2, \theta_1 - \theta_2 \}.
$$
So we have a fundamental cell $C$ of $\frak a'$:
$$
C = \left\{ Y=
\left( \begin{array}{cc} O & X \\ -^tX & O \end{array} \right) \ \bigg| \
X = \left( \begin{array}{cc} \theta_1 & 0 \\ 0 & \theta_2 \end{array} \right),
\begin{array}{c} 0 \leq \theta_1 + \theta_2 \leq \pi \\
0 \leq \theta_1 - \theta_2 \leq \pi \end{array} \right \}.
$$
Thus the isotropy invariants of this case is given by 
$\theta_1 + \theta_2$ and  $\theta_1 - \theta_2$.
It is easy to see that the geometric meaning of $\theta_1 - \theta_2$
is the K\"ahler angle of 2-dimensional subspace $\mbox{Exp} Y$ of $T_o(G/K)$.
On the other hand, there is the other complex structure $J'$
which is defined by
$$
J'e_1 = e_2,\ J'e_2 = - e_1,\ J'e_3 = -e_4,\ J'e_4 = e_3
$$
on $T_o(G/K)$.
We can also check that $\theta_1 + \theta_2$ is the K\"ahler angle
of $\mbox{Exp} Y$ with respect to $J'$.

We attempt to obtain the explicit expression of the Poincar\'e formula
applying the isotropy invariants which we defined above
to Theorem 5.
Let $N$ and $L$ be surfaces in $S^2 \times S^2$.
We take orthonormal bases $\{ u_1, u_2 \}$ and $\{ v_1, v_2 \}$ of
$(dg_x)_o^{-1} (T_x^\perp N)$ and $(dg_y)_o^{-1} (T_y^\perp L)$, respectively.
By the definition, we have
$$
\sigma_K (T_x^\perp N, T_y^\perp L)
= \int_K \| u_1 \wedge u_2 \wedge k^{-1}(v_1 \wedge v_2) \| d\mu_K^{}(k).
$$
Furthermore, by the Hodge star operator,
$$
\sigma_K (T_x^\perp N, T_y^\perp L)
= \int_K | \langle u'_1 \wedge u'_2, k^{-1}(v_1 \wedge v_2) \rangle |
d\mu_K^{}(k),
$$
where $\{ u'_1, u'_2 \}$ is an orthonormal basis of $(dg_x)_o^{-1} (T_x N)$.
We put
$$
a = \left[ \begin{array}{cccc} \cos \phi & \sin \phi & & \\
-\sin \phi & \cos \phi & & \\ & & 1 & \\ & & & 1 \end{array} \right],
\qquad
b = \left[ \begin{array}{cccc} 1 & & & \\ & 1 & & \\
 & & \cos \psi & \sin \psi \\ & & -\sin \psi & \cos \psi \end{array} \right]
$$
and $k = b^{-1}a$, then we have
\begin{equation} \label{eq:3-1}
\sigma_K (T_x^\perp N, T_y^\perp L)
= \int_0^{2 \pi} \int_0^{2 \pi}
| \langle a(u'_1 \wedge u'_2), b(v_1 \wedge v_2) \rangle | d\phi d\psi.
\end{equation}
Since without loss of generalities we can assume that
$(dg_x)_o^{-1} (T_x^\perp N)$ and $(dg_y)_o^{-1} (T_y^\perp L)$ are in
$\mbox{Exp} C$,
we can take $\{ u_1', u_2' \}$ and $\{ v_1, v_2 \}$ as follows:
\begin{eqnarray}
u'_1 = \sin \theta_1 e_1 + \cos \theta_1 e_3, &&
u'_2 = \sin \theta_2 e_2 + \cos \theta_2 e_4, \label{eq:3-2}\\
v_1 = \cos \tau_1 e_1 - \sin \tau_1 e_3, &&
v_2 = \cos \tau_2 e_2 - \sin \tau_2 e_4, \label{eq:3-3}
\end{eqnarray}
with isotropy invariants $\theta_1 \pm \theta_2$ and $\tau_1 \pm \tau_2$.
So we can express the integration of (\ref{eq:3-1})
using $\theta_1$, $\theta_2$, $\tau_1$ and $\tau_2$.
It is complicated to express this general form,
so we shall show for some special cases which is needed
to prove our main theorem.
\begin{thm}
Let $N$ and $L$ be Lagrangian surfaces in $S^2 \times S^2$.
We assume that $L$ is a product of curves in $S^2$.
Then we have
$$
\int_G \sharp (N \cap gL)d\mu(g) = 4 \mathrm{vol}(L)
\int_N \mathrm{length}(\mathrm{Ellip}(\sin^2 \theta_x, \cos^2 \theta_x))
d\mu(x),
$$
where $2\theta_x -\pi/2$ is the K\"ahler angle of $T_x^\perp N$
with respect to $J'$ and
$\mathrm{Ellip}(\alpha, \beta)$ denotes an ellipse defined by
$$
\frac{x^2}{\alpha^2} + \frac{y^2}{\beta^2} = 1.
$$
\end{thm}
\noindent
{\em Proof.}
Since $N$ is a Lagrangian surface, $\theta_1 - \theta_2 = \pi/2$
in (\ref{eq:3-2}), so we put
$$
\theta_1 = \theta, \quad \theta_2 = \theta - \frac{\pi}{2}.
$$
On the other hand,
$L$ is Lagrangian with respect to both $J$ and $J'$,
that is, $\tau_1=\pi/2$ and $\tau_2=0$ in (\ref{eq:3-3}).
Therefore from (\ref{eq:3-1}) we have
$$
\sigma_K (T_x^\perp N, T_y^\perp L)
= \int_0^{2 \pi} \int_0^{2 \pi}
| \cos \phi \cos \psi \cos^2 \theta + \sin \phi \sin \psi \sin^2 \theta |
d\phi d\psi.
$$
In~\cite{Kang}, Kang calculated this type of integrals directly
and expressed it by elliptic functions.
But we give here a geometrical simple computation.
Now we take a subspace which is given by
$$
V = \mbox{span} \{ e_2 \wedge e_3, e_2 \wedge e_4 \}
$$
in $\bigwedge^2 \big( T_o(G/K) \big)$.
Then $b(v_1 \wedge v_2)$ moves on the unit circle in $V$
with the parameter $\psi$.
Let $P$ denote the orthogonal projection
from $\bigwedge^2 \big( T_o(G/K) \big)$ to $V$.
>From (\ref{eq:3-1}) we have
$$
\sigma_K (T_x^\perp N, T_y^\perp L)
= \int_0^{2 \pi} \int_0^{2 \pi}
| \langle P(a(u'_1 \wedge u'_2)), b(v_1 \wedge v_2) \rangle | d\phi d\psi.
$$
Here $P(a(u'_1 \wedge u'_2))$ moves  with parameter $\phi$
on the ellipse which defined by
$$
\frac{x^2}{\cos^4 \theta} + \frac{y^2}{\sin^4 \theta} = 1
$$
in $V$.
Hence we put
$$
r_\phi = \sqrt{\cos^2 \phi \cos^4 \theta + \sin^2 \phi \sin^4 \theta},
$$
then we have
\begin{eqnarray*}
\sigma_K (T_x^\perp N, T_y^\perp L)
&=& \int_0^{2 \pi} \int_0^{2 \pi} | r_\phi \cos \psi | d\psi d\phi \\
&=& \int_0^{2 \pi} r_\phi d\phi \int_0^{2 \pi} | \cos \psi | d\psi \\
&=& 4 \mbox{length}(\mbox{Ellip}(\sin^2 \theta, \cos^2 \theta)).
\end{eqnarray*}
Thus we complete the proof from Theorem 5. \hfill \qed

\vspace*{12pt}
Proposition 4 is immediately obtained from Theorem 6.

\section{Proof of the main theorem}

\noindent
{\em Proof of Theorem 2.}
Let $L:=S^1 \times S^1$ be a totally geodesic Lagrangian torus in
$S^2 \times S^2$.
Let $\rho$ be a Hamiltonian isotopy of $S^2 \times S^2$.
By inequalities (1) and (2), we have
\begin{eqnarray*}
16 \mathrm{vol}(\rho(L))\mathrm{vol}(L)
&\geq& \int_{SO(3) \times SO(3)} \sharp (\rho(L) \cap gL)d\mu(g) \\
&=& \int_{SO(3) \times SO(3)} \sharp (g^{-1} \circ \rho(L) \cap L)d\mu(g) \\
&\geq& \int_{SO(3) \times SO(3)}
 \sum_{i=0}^{2} \mathrm{rank} H_i(L, {\mathbb Z}/2{\mathbb Z})d\mu(g) \\
&=& 4 \mathrm{vol}(SO(3) \times SO(3)).
\end{eqnarray*}
Since $\mathrm{vol}(SO(3))=8\pi^2$ and $\mathrm{vol}(L)=4\pi^2$, we have
\begin{eqnarray*}
\mathrm{vol}(\rho(L)) \geq 4\pi^2=\mathrm{vol}(L).
\end{eqnarray*}
\hfill \qed

\section*{Acknowledgements}

We would like to thank Professor Yong-Guen Oh for some helpful
comments on the original version of this paper.

\noindent
Hiroshi Iriyeh\\
{\sc Department of Mathematics, Graduate School of Science,\\
Tokyo Metropolitan University, Minami-Ohsawa 1-1,
Hachioji, Tokyo 192-0397, Japan}

{\it E-mail address} : {\tt hirie@comp.metro-u.ac.jp}

\vspace{5mm}

\noindent
Hajime Ono\\
{\sc Department of Mathematics, Graduate School of Science,\\
Tokyo Metropolitan University, Minami-Ohsawa 1-1,
Hachioji, Tokyo 192-0397, Japan}

{\it E-mail address} : {\tt onola@comp.metro-u.ac.jp}

\vspace{5mm}

\noindent
Takashi Sakai\\
{\sc Department of Mathematics, Graduate School of Science,\\
Tokyo Metropolitan University, Minami-Ohsawa 1-1,
Hachioji, Tokyo 192-0397, Japan}

{\it E-mail address} : {\tt tsakai@comp.metro-u.ac.jp}

\end{document}